\documentclass[pdflatex,sn-mathphys-num]{sn-jnl}

\usepackage{graphicx}%
\usepackage{multirow}%
\usepackage{amsmath,amssymb,amsfonts}%
\usepackage{amsthm}%
\usepackage{mathrsfs}%
\usepackage[title]{appendix}%
\usepackage{xcolor}%
\usepackage{textcomp}%
\usepackage{manyfoot}%
\usepackage{booktabs}%
\usepackage{algorithm}%
\usepackage{algorithmicx}%
\usepackage{algpseudocode}%
\usepackage{listings}%
\usepackage[T1]{fontenc}
\usepackage[utf8]{inputenc}

\theoremstyle{thmstyleone}%
\theoremstyle{thmstyletwo}%

\theoremstyle{thmstylethree}%

\newtheorem{theo}{Theorem}
\newtheorem{lem}[theo]{Lemma}
\newtheorem{prop}[theo]{Proposition}
\newtheorem{cor}[theo]{Corollary}
\newtheorem{rem}{Remark}

\newenvironment{dem}[1][Proof]{\noindent \textbf{#1.} }{\ \rule{0.5em}{0.5em}}
\begin{document}

\title[Normal cones to sublevel sets]{Normal cones to sublevel sets of convex and quasi-convex supremum
functions}


\author[1]{\fnm{Stephanie} \sur{Caro}}\email{stcaro@unap.cl}
\equalcont{These authors contributed equally to this work.}

\author[2]{\fnm{Rafael} \sur{Correa}}\email{rcorrea@dim.uchile.cl}
\equalcont{These authors contributed equally to this work.}

\author*[3]{\fnm{Abderrahim} \sur{Hantoute}}\email{hantoute@ua.es}
\equalcont{These authors contributed equally to this work.}

\affil[1]{\orgdiv{Faculty of Sciences}, \orgname{Institute of Exact and Natural Sciences, University Arturo Prat}, \orgaddress{
\city{Iquique}, 
\country{Chile}}}

\affil[2]{\orgdiv{Mathematical Engineering Department}, \orgname{University of Chile}, \orgaddress{
\city{Santiago}, 
\country{Chile}}}

\affil*[3]{\orgdiv{Mathematical Department}, \orgname{University of Alicante}, \orgaddress{
\city{Alicante}, 
\country{Spain}}}



\date{\today}
\maketitle

\begin{abstract}

We provide sharp and explicit characterizations of the normal cone to sublevel
sets of suprema of arbitrary functions, expressed exclusively in terms of
subdifferentials of the data functions. In the convex case, the resulting
formulas involve the approximate subdifferential of the individual data
functions at the nominal point. In contrast, the quasi-convex framework
requires the use of the Fréchet subdifferential of these data functions but
evaluated at nearby points.\ These results are applied to derive optimality
conditions for infinite convex and quasi-convex optimization problems.



\end{abstract}

\section{Introduction}

The normal cone to sublevel sets plays a central role in optimization and
variational analysis, both from a theoretical and an applied perspective. As
documented in many classical references (e.g., \cite{DiGoLo06, HiMoVoSe39,
Mor, Ro70, RoWe} and references therein), its importance stems from the fact
that a wide range of constraint systems, optimality conditions of first and
higher orders, duality frameworks, as well as fundamental concepts and
techniques of variational and nonsmooth analysis, and stability properties can
be formulated and analyzed entirely in terms of normal cones to sets of the
form
\[
\lbrack f\leq0]:=\{x\in X:f(x)\leq0\},
\]
where $f:X\rightarrow\overline{\mathbb{R}}:=\mathbb{R\cup\{\pm\infty\}}$ is a
suitable function defined on a linear\ space $X,$ encompassing the data of the
underlying optimization or variational problem. For instance, in infinite
optimization problems involving infinitely many constraints, first-order
necessary optimality conditions are naturally expressed using normal cones. If
$x$ is a local minimizer of a function $f_{0}:X\rightarrow\overline
{\mathbb{R}}$ over a constraint set $C:=\{x\in X:f_{t}(x)\leq0,$ $t\in T\},$
where $T$ is an arbitrary index set, then under suitable qualification
conditions one typically has%
\[
0\in\partial f_{0}(x)+\mathrm{N}_{C}(x),
\]
where $\mathrm{N}_{C}(x)$ is the normal cone to $C,$ and $\partial f_{0}$
stands for\ an appropriate subdifferential. In this setting, it is often
convenient to represent the constraints system by the pointwise supremum
function (see, e.g., \cite{CHLBook, HLZ08}):%
\[
f:=\sup_{t\in T}f_{t}.
\]
In this context, describing\ the normal cone to $[f\leq0]$ through explicit
representations in terms of subdifferentials of the data functions $f_{t}$ is
essential for deriving optimality conditions, sensitivity results, and dual
formulations. Different characterizations of this normal cone have been
proposed\ in recent years. In particular, \cite{CaTh14} provides limiting
sequential-type characterizations in the case of single functions defined on
Banach spaces.\ Extensions to functions expressed as suprema are studied in
\cite{HaLo23, HaSv17} (see, also, \cite{CHLBook}), where the characterizations
involve subdifferentials of convex combinations of the data functions. 

In this paper, we establish new characterizations of the normal cone to the
sublevel set of arbitrary pointwise suprema. In the convex setting, the
resulting\ representations are expressed exclusively in terms of approximate
subdifferentials of the\ individual data functions evaluated at the nominal
points.\ In contrast to \cite{HaLo23, HaSv17}, our characterization involve
the data functions themselves rather than convex combinations thereof. In a
second step, we extend these results to the quasi-convex setting, aiming to
better understand the geometric structure of the normal cones to sublevel sets
of suprema of quasi-convex functions, which play a fundamental role in
quasi-convex analysis (see, e.g., \cite{Da} and references therein). Two types
of characterizations are obtained in the quasi-convex framework: a global one,
expressed in terms of normal cones to the sublevel sets of the individual data
functions, and a local one, formulated via the Fréchet subdifferential of
these data functions evaluated at nearby points. These explicit descriptions
of the normal cone therefore will lead us directly to computable and
verifiable optimality conditions, particularly for optimization problems
involving\ infinitely many convex and quasi-convex constraints. Our approach
relies on techniques from convex and nonsmooth analysis, notably a recent
characterization of the normal cone to the effective domain of pointwise
suprema established\ in\ \cite{CaHa26}.

This paper is organized as follows: Section \ref{sec2} introduces the main
notation and recalls the preliminary results that will be used in the sequel.
In Section \ref{sec4}, we establish the characterizations of the normal cone
to sublevel sets of convex functions and derive optimality conditions for
infinite convex optimization. Section \ref{sec5} is devoted to extensions of
these results to the quasi-convex setting. Finally, concluding remarks are
presented\ in Section \ref{sec6}.

\section{Notation and preliminaries\label{sec2}}

We consider\ a real locally convex space $X$ (lcs, for short), whose
topological dual is denoted by $X^{\ast}.$ The duality$\ $pairing is given by
\[
(x^{\ast},x)\in X^{\ast}\times X\mapsto\langle x^{\ast},x\rangle:=x^{\ast}(x),
\]
and both $X$ and $X^{\ast}$ are endowed with any topology compatible with this
duality.$\ $The zero vector in both $X$ and $X^{\ast}$ is denoted by $\theta.$
The family of all closed, convex, and balanced neighborhoods of\textbf{\ }%
$\theta\ $(called\textbf{\ }$\theta$\textbf{-}neighborhoods\textbf{)} is
denoted\ by\textbf{\ }$\mathcal{N}$.

We use the notation $\overline{\mathbb{R}}:=\mathbb{R}\cup\{-\infty,+\infty\}$
and $\mathbb{R}_{\infty}:=\mathbb{R}\cup\{+\infty\}$, and adopt the
conventions\emph{\ }$\left(  +\infty\right)  +(-\infty)=+\infty$ and
$0.(+\infty):=+\infty$. 

For subsets\ $A$, $B\subset X,$ their Minkowski sum is defined by
$A+B:=\{a+b:\ a\in A,\text{ }b\in B\}$, with the convention $A+\emptyset
=\emptyset+A=\emptyset$. The convex hull and closure of a set $A$ are denoted
by $\operatorname*{co}(A)$ and $\operatorname{cl}(A)$ (or $\bar{A}$),
respectively. If $X=\mathbb{R}^{n},$ we denote by $\operatorname*{ri}(A)$ the
relative interior of $A$. The\ recession cone of\ a nonempty closed convex set
$A$\ is
\[
\lbrack A]_{\infty}:=\cap_{t>0}t(A-a),
\]
for an arbitrary $a\in A;$ by convention, we set $[\emptyset]_{\infty
}=\{\theta\}.$

Given a function $f:X\rightarrow\overline{\mathbb{R}}$, we denote by
$\operatorname*{dom}f$ and $\operatorname*{epi}f$ its (effective) domain and
epigraph,$\ $respectively. We say that $f$ is proper if $\operatorname*{dom}%
f\neq\emptyset$ and $f(x)>-\infty,$ for all $x\in X$. For\ $c\in\mathbb{R}$,
we define the sublevel sets
\[
\lbrack f\leq c]:=\{x\in X:f(x)\leq c\},
\]
and $[f<c]:=\{x\in X:f(x)<c\}.$ The function $f$ is said to be convex if
$\operatorname*{epi}f$ is a convex set in $X\times\mathbb{R}$, and lower
semicontinuous (lsc, for short) if $\operatorname*{epi}f$ is closed. It is
called\ quasi-convex if $[f\leq c]$ is convex, for every $c\in\mathbb{R}.$ The
closed hull of\ $f,$ denoted by $\operatorname*{cl}f$ (or $\overline{f}$), is
the function whose epigraph is $\operatorname*{cl}(\operatorname*{epi}f).$

For a set $A\subset X^{\ast}$ (or $X$), its support function is defined by
\[
\sigma_{A}:=\sup_{a\in A}\left\langle a,\cdot\right\rangle .
\]
Given a set $U,$ the indicator function of a set $A\subset U,$ denoted by
$\mathrm{I}_{A}:U\rightarrow\mathbb{R}_{\infty}$, is defined by $\mathrm{I}%
_{A}(a):=0$ if $a\in A,$ and $\mathrm{I}_{A}(a):=+\infty$ if $a\in U\setminus
A$. We adopt\ the conventions $\sigma_{\emptyset}\equiv-\infty$ and
$\mathrm{I}_{\emptyset}\equiv+\infty.$

If $A\subset X$ is a nonempty set, then (see, e.g., \cite[Eq. (3.52), p.
80]{CHLBook})
\begin{equation}
\operatorname{cl}(\operatorname{dom}\sigma_{A})=([\overline{\operatorname*{co}%
}\left(  A\right)  ]_{\infty})^{-}.\label{rl}%
\end{equation}
Since $\operatorname{dom}\sigma_{A}=\operatorname{dom}\sigma_{A\cup\{\theta
\}}$, the last\ relation together with the bipolar theorem, implies
\[
\lbrack\overline{\operatorname*{co}}(A)]_{\infty}=[\overline
{\operatorname*{co}}(A\cup\{\theta\})]_{\infty}.
\]
Hence, when dealing with the recession cone $[\overline{\operatorname*{co}%
}(A)]_{\infty}$, we may always assume without loss of generality that
$\theta\in A.$ Consequently, for any family of nonempty closed convex sets
$A_{i}\subset X,$ $i\in I,$ it follows that\
\[
\lbrack A_{i}]_{\infty}=[\overline{\operatorname*{co}}(A_{i}\cup
\{\theta\})]_{\infty}\subset\overline{\operatorname*{co}}(A_{i}\cup
\{\theta\}).
\]
Taking into account \cite[Eq. (2.23), p. 25]{CHLBook} and the convention
$[\emptyset]_{\infty}=\{\theta\},$ we obtain
\begin{equation}
\cap_{i\in I}[A_{i}]_{\infty}=\cap_{i\in I}[A_{i}\cup\{\theta\}]_{\infty
}=[\cap_{i\in I}(A_{i}\cup\{\theta\})]_{\infty}=[\cap_{i\in I}A_{i}]_{\infty
}.\label{ricone}%
\end{equation}
Given $\varepsilon\geq0$, the $\varepsilon$-subdifferential of $f:X\rightarrow
\overline{\mathbb{R}}$ at a point $x\in X$ is defined by
\[
\partial_{\varepsilon}f(x):=\{x^{\ast}\in X^{\ast}:\langle x^{\ast}%
,y-x\rangle\leq f(y)-f(x)+\varepsilon,\text{ for all }y\in X\}
\]
whenever $x\in f^{-1}(\mathbb{R}),$ and $\partial_{\varepsilon}f(x):=\emptyset
$ otherwise. The (exact) subdifferential of $f$\ at $x$ is given by
\[
\partial f(x):=\partial_{0}f(x)=\cap_{\varepsilon>0}\partial_{\varepsilon
}f(x).
\]
We call $\mathrm{N}_{A}^{\varepsilon}(x):=\partial_{\varepsilon}\mathrm{I}%
_{A}(x)$\ the $\varepsilon$-normal set to $A\subset X$ at $x\in X$; in
particular, $\mathrm{N}_{A}(x):=\mathrm{N}_{A}^{0}(x)$ is the normal cone\ to
$A\subset X$ at $x$. Note that\ $\mathrm{N}_{A}^{\varepsilon}(x)=\mathrm{N}%
_{\bar{A}}^{\varepsilon}(x)$ whenever $x\in A.$ Consequently,
since\ $\operatorname*{cl}(\operatorname*{dom}f)=\operatorname*{cl}%
(\operatorname*{dom}\bar{f})$ for every function $f$\ on $X$, it follows
that\
\[
\mathrm{N}_{\operatorname*{dom}f}^{\varepsilon}(x)=\mathrm{N}%
_{\operatorname*{dom}\bar{f}}^{\varepsilon}(x),\text{ \ for all }%
x\in\operatorname*{dom}f\text{ and }\varepsilon\geq0.
\]
If $(X,\left\Vert \cdot\right\Vert )$ is\ a Banach space and $f:X\rightarrow
\mathbb{R}_{\infty},$ an element $x^{\ast}\in X^{\ast}$ is called a
Fréchet-subgradient of $f$ at $x\in\operatorname*{dom}f$ if
\[
\liminf_{y\rightarrow x}\frac{f(y)-f(x)-\left\langle x^{\ast},y-x\right\rangle
}{\left\Vert y-x\right\Vert }\geq0.
\]
The set of all such elements is called the Fréchet-subdifferential of $f$ at
$x.$ We denote it by $\partial f(x),$ as in the convex setting, because it
reduces to the Fenchel subdifferential whenever $f$ is convex.

Next, we recall some topological properties of sublevel sets associated
with\ a convex function $f:X\rightarrow\overline{\mathbb{R}}$. If $[f<0]\cap
f^{-1}(\mathbb{R})\neq\emptyset,$ then (e.g., \cite[Lemma 3.3.3]{CHLBook})
\begin{equation}
\operatorname*{cl}([f<0])=\operatorname*{cl}([f\leq0])=[\bar{f}\leq
0].\label{lemnewa}%
\end{equation}
Moreover, we always have that
\begin{equation}
\operatorname*{cl}([f\leq0])=[\bar{f}\leq0],\label{lemnewaa}%
\end{equation}
when either $X=\mathbb{R}^{n}$ and $\operatorname*{ri}(\operatorname*{dom}%
f)\cap\lbrack f\leq0]\neq\emptyset,$ or $f$ is continuous at some point in
$[f\leq0].$ Indeed, consider the convex functions $f_{s}:=sf,$ $s>0,$ so that
\[
\mathrm{I}_{[f\leq0]}(x)=\sup_{s>0}f_{s}(x),\text{ for all }x\in X.
\]
Since $\operatorname*{dom}(\sup_{s>0}f_{s})=[f\leq0]\neq\emptyset,$ and
$\operatorname*{cl}(\sup_{s>0}f_{s})=\sup_{s>0}\overline{f_{s}}\ $by
\cite[Proposition 5.2.4((i) and (iii))]{CHLBook} , we have
\[
\mathrm{I}_{\operatorname*{cl}([f\leq0])}=\operatorname*{cl}(\mathrm{I}%
_{[f\leq0]})=\operatorname*{cl}(\sup_{s>0}f_{s})=\sup_{s>0}\overline{f_{s}%
}=\sup_{s>0}s\bar{f}=\mathrm{I}_{[\bar{f}\leq0]},
\]
which yields (\ref{lemnewaa}).

The following lemma is the counterpart of (\ref{lemnewaa}) for quasi-convex
functions. This result was given in \cite[Proposition 3.4 ]{Cr2005} under the
condition\ $\operatorname*{int}([f\leq0])\not =\emptyset.$

\begin{lem}
\label{lemnewb}Let $f:X\rightarrow\overline{\mathbb{R}}$ be a quasi-convex
function such that $[f\leq0]\neq\emptyset.$\ If either $X=\mathbb{R}^{n}$ and
$\operatorname*{ri}(\operatorname*{dom}f)\cap\lbrack f\leq0]\neq\emptyset$, or
$f$ is continuous at some point in $[f\leq0],$ then\
\[
\operatorname*{cl}([f\leq0])=[\bar{f}\leq0].
\]

\end{lem}

\begin{dem}
Consider the functions $g_{s}:=\mathrm{I}_{[-\infty,s]}\circ f,$ $s>0.$ Then
each $g_{s}$ is convex and we have that
\[
\mathrm{I}_{[f\leq0]}(x)=\sup_{s>0}g_{s}(x),\text{ for all }x\in X.
\]
Assume that\ $f$ is continuous at some point $x_{0}\in\lbrack f\leq0].$
Then\ each $g_{s},$ $s>0,$ is also continuous at $x_{0}.$ Hence,
\cite[Proposition 5.2.4(i)]{CHLBook}\ entails
\[
\mathrm{I}_{\operatorname*{cl}([f\leq0])}=\operatorname*{cl}(\mathrm{I}%
_{[f\leq0]})=\operatorname*{cl}(\sup_{s>0}g_{s})=\sup_{s>0}(\operatorname*{cl}%
g_{s})=\sup_{s>0}(\operatorname*{cl}(\mathrm{I}_{[-\infty,s]}\circ f)).
\]
Therefore, since $[\bar{f}\leq0]\subset\operatorname*{cl}([f\leq s])$ for all
$s>0,$ we obtain
\[
\mathrm{I}_{\operatorname*{cl}([f\leq0])}\leq\mathrm{I}_{[-\infty,0]}\circ
\bar{f}=\mathrm{I}_{[\bar{f}\leq0]};
\]
that is, $[\bar{f}\leq0]\subset\operatorname*{cl}([f\leq0]).$ As the opposite
inclusion is straightforward, the conclusion follows under the given
continuity assumption on $f.$

Assume now that $X=\mathbb{R}^{n}$ and $\operatorname*{ri}(\operatorname*{dom}%
f)\cap\lbrack f\leq0]\neq\emptyset.$ We pick any point $z_{0}\in\lbrack
f\leq0]$ and consider the convex function $\tilde{f}:E:=\operatorname*{aff}%
(\operatorname*{dom}f-z_{0})\rightarrow\overline{\mathbb{R}}$ defined by
$\tilde{f}(z):=f(z-z_{0}),$ $z\in E.$ Then $\operatorname*{cl}(\tilde
{f})(z):=\bar{f}(z-z_{0})$ for all $z\in E$, and we check that $[\tilde{f}%
\leq0]=z_{0}+[f\leq0]$ and $[\operatorname*{cl}(\tilde{f})\leq0]=z_{0}%
+[\bar{f}\leq0].$ Consequently, since $\tilde{f}$ is continuous at every point
in $z_{0}+\operatorname*{ri}(\operatorname*{dom}f),$ the first paragraph
entails that
\[
\lbrack\bar{f}\leq0]=[\operatorname*{cl}(\tilde{f})\leq0]-z_{0}%
=\operatorname*{cl}([\tilde{f}\leq0])-z_{0}=\operatorname*{cl}([f\leq0]),
\]
which is the desired relation.
\end{dem}

The following corollary is specific to the one-dimensional case and does not
hold in higher dimensions, even for convex functions.

\begin{cor}
Let\ $X=\mathbb{R}$. Then any\ quasi-convex function $f:\mathbb{R}%
\rightarrow\overline{\mathbb{R}}$ with $[f\leq0]\neq\emptyset\ $satisfies
\[
\operatorname*{cl}([f\leq0])=[\bar{f}\leq0].
\]

\end{cor}

\begin{dem}
If $\operatorname*{ri}(\operatorname*{dom}f)\cap\lbrack f\leq0]\neq\emptyset,$
then the conclusion follows by Lemma \ref{lemnewb}. We so assume that
$\operatorname*{ri}(\operatorname*{dom}f)\cap\lbrack f\leq0]=\emptyset.$ Since
both $\operatorname*{ri}(\operatorname*{dom}f)$ and $[f\leq0]$ are nonempty
intervals of $\mathbb{R}$, it follows that $[f\leq0]$ reduces to\ a singleton,
say $\bar{x}.$ Assume, by contradiction, that $[\bar{f}\leq0]$ contains
another point distinct from $\bar{x},$ say $\tilde{x}.$ Then, since $\bar{f}$
is also quasi-convex (\cite[p. 128 ]{Cr2005}), $\bar{f}$ is continuous almost
every on the interval $[\bar{x},\tilde{x}].$ Thus, by Proposition
\cite[Proposition 3.5 ]{Cr2005}, the function $f$ is also continuous almost
every on the interval $[\bar{x},\tilde{x}].$ In particular, there exists a
point\ $\hat{x}\in]\bar{x},\tilde{x}[$ where $f$ is continuous. So, $f(\hat
{x})=\bar{f}(\hat{x})\leq0$ and we get $\hat{x}\in\lbrack f\leq0]\setminus
\{\bar{x}\},$ contradicting the fact that $[f\leq0]=\{\bar{x}\}.$ Therefore,
$[\bar{f}\leq0]$ is a singleton, and the trivial inclusion $\{\bar{x}%
\}=[f\leq0]\subset\lbrack\bar{f}\leq0]$ implies that $\operatorname*{cl}%
([f\leq0])=[\bar{f}\leq0].\medskip$
\end{dem}

We close this section of preliminaries by recalling a characterization from
\cite{CaHa26} of the normal cone $\mathrm{N}_{\operatorname*{dom}f}(x)$ for
$x\in\operatorname*{dom}f,$ where $f$ is the pointwise supremum of convex
functions $f_{t}:X\rightarrow\mathbb{R}_{\infty},$ $t\in T.$
This\ characterization relies on the $\varepsilon$-subdifferential of the
functions that are $\varepsilon$-active at $x;$ i.e., those indexed in
\begin{equation}
T_{\varepsilon}(x):=\{t\in T:f_{t}(x)\geq f(x)-\varepsilon\},\text{
}\varepsilon>0,\label{teps2}%
\end{equation}
as well as on the $\varepsilon$-subdifferential of the\ non-$\varepsilon
$-active functions, which enter with appropriate weights. Nonproper functions
also play a role in this characterization through their contributions via the
$\varepsilon$-normal sets $\mathrm{N}_{\operatorname*{dom}f_{t}}^{\varepsilon
}(x).$

\begin{prop}
\label{Theo13} Given lsc convex\ functions $f_{t}:X\rightarrow\overline
{\mathbb{R}},\ t\in T,$\ $f:=\sup_{t\in T}f_{t},$ and $x\in\operatorname*{dom}%
f,$ we denote
\[
\mathcal{P}:=\{t\in T:f_{t}\text{ is proper}\}
\]
and
\begin{equation}
\rho_{t,\varepsilon,x}:=\frac{-\varepsilon}{2f_{t}(x)-2f(x)+\varepsilon}\text{
if }t\in T\setminus T_{\varepsilon}(x),\text{ }\rho_{t,\varepsilon,x}:=1\text{
if\ }t\in T_{\varepsilon}(x).\label{vetb}%
\end{equation}
Then we have
\[
\mathrm{N}_{\operatorname*{dom}f}(x)=\left[  \overline{\operatorname*{co}%
}\left(  \left(  {%
{\textstyle\bigcup\limits_{t\in\mathcal{P}}}
}\partial_{\varepsilon}(\alpha_{t}f_{t})(x)\right)  {\cup}\left(  {%
{\textstyle\bigcup\limits_{t\in T\backslash\mathcal{P}}}
}\mathrm{N}_{\operatorname*{dom}f_{t}}^{\varepsilon}(x)\right)  \right)
\right]  _{\infty},
\]
for any\ $\varepsilon>0$ and $\alpha_{t}\in\mathbb{R}_{+}$ such that
\[
\alpha_{t}\geq\rho_{t,\varepsilon,x},\text{ for all }t\in\mathcal{P},
\]%
\[
\sup_{t\in\mathcal{P}}\alpha_{t}<\infty,\text{ and }\inf_{t\in\mathcal{P}%
}\alpha_{t}f_{t}(x)>-\infty.
\]

\end{prop}

\section{Normal cone to sublevel sets in the convex case\label{sec4}}

We consider an arbitrary family of convex functions $f_{t}:X\rightarrow
\overline{\mathbb{R}}$, $t\in T$, and the associated supremum function
\[
f:=\sup_{t\in T}f_{t}.
\]
Our first objective is to characterize the normal cone $\mathrm{N}_{[f\leq
0]}(x)$ at points $x\in\lbrack f\leq0]$ exclusively in terms of the data
functions $f_{t},$ through their $\varepsilon$-subdifferentials. These results
will provide the foundation for deriving new optimality conditions in infinite
convex optimization.

We state below the first main result of this section, which relies on the
following lower semicontinuity-type\ condition
\begin{equation}
\overline{\lbrack f\leq0]}=\cap_{t\in T}[\bar{f}_{t}\leq0],\label{cc2}%
\end{equation}
involving the lsc hulls of $f_{t}$ (see Remark \ref{remm}$(ii)$ for a
discussion\ of this assumption). We also\ use the notation\
\begin{equation}
T_{\varepsilon/s}(x):=\{t\in\mathcal{P}:sf_{t}(x)\geq-\varepsilon\},\text{
}\varepsilon,s>0,\label{az}%
\end{equation}
where
\[
\mathcal{P}:=\{t\in T:\bar{f}_{t}\text{ proper}\}.
\]

\begin{theo}
\label{thmcone} Let\ $f_{t}:X\rightarrow\overline{\mathbb{R}},$ $t\in T,$ be
convex, and denote $f:=\sup_{t\in T}f_{t}.$ Assume that condition (\ref{cc2})
holds. Then, for every $x\in\lbrack f\leq0]$ such that $f(x)\in\mathbb{R},$ we
have
\begin{equation}
\mathrm{N}_{[f\leq0]}(x)=\left[  \overline{\operatorname*{co}}\left(  \left(
{\textstyle\bigcup\limits_{t\in T_{\varepsilon/s}(x),s>0\text{ }}}
\partial_{\varepsilon}(sf_{t})(x)\right)  \cup\left(
{\textstyle\bigcup\limits_{t\in T\setminus\mathcal{P}\text{ }}}
\mathrm{N}_{\operatorname*{dom}f_{t}}^{\varepsilon}(x)\right)  \right)
\right]  _{\infty},\text{ for all }\varepsilon>0.\label{a1a}%
\end{equation}
Consequently,
\begin{equation}
\mathrm{N}_{[f\leq0]}(x)=%
{\textstyle\bigcap\limits_{\varepsilon>0}}
\overline{\operatorname*{co}}\left(  \left(
{\textstyle\bigcup\limits_{t\in T_{\varepsilon/s}(x),s>0\text{ }}}
\partial_{\varepsilon}(sf_{t})(x)\right)  \cup\left(
{\textstyle\bigcup\limits_{t\in T\setminus\mathcal{P}\text{ }}}
\mathrm{N}_{\operatorname*{dom}f_{t}}^{\varepsilon}(x)\right)  \right)
.\label{a2b}%
\end{equation}

\end{theo}

\begin{dem}
Fix\ $x\in(\operatorname*{dom}f)\cap\lbrack f\leq0]$ and $\varepsilon>0.$ We
first suppose, in addition to the current assumptions,\ that all the $f_{t}$'s
are lsc. Then the functions
\[
f_{t,s}:=sf_{t},\text{ }(t,s)\in T\times]0,+\infty\lbrack,
\]
are convex and lsc, and we have
\[
g:=\sup_{t\in T,s>0}f_{t,s}=\mathrm{I}_{[f\leq0]}=\mathrm{I}%
_{\operatorname*{cl}([f\leq0])}.
\]
It follows that\ $g(x)=0$ and $\operatorname*{dom}g=[f\leq0].$ Moreover, the
inequality $f_{t,s}(x)\geq g(x)-\varepsilon$ holds if and only if $t\in
T_{\varepsilon/s}(x)$ with $s>0.$ Let us\ consider the weights\ $\rho
_{t,s}:=\rho_{t,s,\varepsilon,x}\ $associated with the family $\{f_{t,s}:$
$t\in T,$ $s>0\}$ (see (\ref{vetb})), defined\ by
\[
\rho_{t,s}:=1\text{ if }t\in T_{\varepsilon/s}(x)\text{ and}\ s>0,
\]
and
\[
\rho_{t,s}:=\frac{-\varepsilon}{2f_{t,s}(x)-2g(x)+\varepsilon}=\frac
{-\varepsilon}{\varepsilon+2sf_{t}(x)},\text{ if }t\in T\setminus
T_{\varepsilon/s}(x)\text{ and }s>0.
\]
We also choose\ the positive numbers $\alpha_{t,s}$ defined by
\[
\alpha_{t,s}:=-\varepsilon/(sf_{t}(x)),\text{ for }t\in\mathcal{P}\setminus
T_{\varepsilon/s}(x)\text{ and }s>0,
\]
and
\[
\alpha_{t,s}:=1,\text{ for}\ t\in T_{\varepsilon/s}(x)\text{ and }s>0.
\]
Observe that $f_{t}(x)<0$ for all $t\in T\setminus T_{\varepsilon/s}(x)$ and
$s>0.$ Then, for $t\in\mathcal{P}\setminus T_{\varepsilon/s}(x)$ and $s>0,$ we
have
\[
\alpha_{t,s}=\varepsilon/(-sf_{t}(x))\leq1,\text{ }%
\]%
\[
\frac{-\varepsilon}{(\varepsilon+2sf_{t}(x))\alpha_{t,s}}=\frac{sf_{t}%
(x)}{\varepsilon+2sf_{t}(x)}<1,\text{ and }\alpha_{t,s}f_{t,s}(x)=-\varepsilon
.
\]
This shows that $\alpha_{t,s}\geq\rho_{t,s}$ (for all $t\in\mathcal{P}$ and
$s>0$), $\sup_{t\in\mathcal{P}\text{, }s>0}\alpha_{t,s}\in\mathbb{R},$ and
\[
\inf_{t\in\mathcal{P}\text{, }s>0}\alpha_{t,s}f_{t,s}(x)=\min\{\inf_{t\in
T_{\varepsilon/s}(x)\text{, }s>0}f_{t,s}(x),\inf_{t\in\mathcal{P}\setminus
T_{\varepsilon/s}(x)\text{, }s>0}(-\varepsilon)\}=-\varepsilon\in\mathbb{R}.
\]
Therefore, applying Proposition\ \ref{Theo13} to the family $\{f_{t,s}%
:=sf_{t},$ $(t,s)\in T\times]0,+\infty\lbrack\}$ yields
\begin{align}
\mathrm{N}_{[f\leq0]}(x)  & =\mathrm{N}_{\operatorname*{dom}g}(x)\nonumber\\
& =\left[  \overline{\operatorname*{co}}\left(  \left(  \cup_{t\in
T_{\frac{\varepsilon}{s}}(x),s>0\text{ }}\partial_{\varepsilon}(\alpha
_{t}f_{t,s})(x)\right)  \cup\left(  \cup_{t\in T\setminus\mathcal{P}\text{,
}s>0\text{ }}\mathrm{N}_{\operatorname*{dom}f_{t,s}}^{\varepsilon}(x)\right)
\right)  \right]  _{\infty}\nonumber\\
& =\left[  \overline{\operatorname*{co}}\left(  A_{\varepsilon}\cup
B_{\varepsilon}\cup\left(  \cup_{t\in T\setminus T_{\varepsilon/s}%
(x),s>0\text{ }}\partial_{\varepsilon}\left(  r_{t}f_{t}\right)  (x)\right)
\right)  \right]  _{\infty},\text{ }\label{yarab}%
\end{align}
where we denoted $r_{t}:=-\varepsilon/f_{t}(x)$ for $t\in T\setminus
T_{\varepsilon/s}(x)$ and $s>0$, and
\begin{equation}
A_{\varepsilon}:=\cup_{t\in T_{\varepsilon/s}(x),s>0\text{ }}\partial
_{\varepsilon}(sf_{t})(x),\text{ }B_{\varepsilon}:=\cup_{t\in T\setminus
\mathcal{P}\text{ }}\mathrm{N}_{\operatorname*{dom}f_{t}}^{\varepsilon
}(x).\label{yarab2}%
\end{equation}
More precisely, since $r_{t}>0$ and $r_{t}f_{t}(x)=-\varepsilon,$ for all
$T\setminus T_{\varepsilon/s}(x)$ and $s>0,$ we have
\[
\cup_{t\in T\setminus T_{\varepsilon/s}(x),s>0\text{ }}\partial_{\varepsilon
}\left(  r_{t}f_{t}\right)  (x)\subset\cup_{t\in T_{\varepsilon/s}%
(x),s>0\text{ }}\partial_{\varepsilon}\left(  sf_{t}\right)
(x)=A_{\varepsilon},
\]
and (\ref{yarab}) simplifies\ to the desired relation (\ref{a1a}) when all the
$f_{t}$'s are lsc.

We now address the general case by considering\ the lsc convex hulls\ $\bar
{f}_{t},$ $t\in T.$ Since $x\in\lbrack f\leq0]\subset\lbrack\bar{f}\leq0]$,
condition (\ref{cc2}) ensures that
\[
\mathrm{N}_{\overline{[f\leq0]}}(x)=\mathrm{N}_{[\sup_{t\in T}\bar{f}_{t}%
\leq0]}(x),
\]
and so, applying the first paragraph with $\varepsilon/2$ to the family
$\{\bar{f}_{t},$ $t\in T\},$\
\begin{equation}
\mathrm{N}_{\overline{[f\leq0]}}(x)=\left[  \overline{\operatorname*{co}%
}\left(  \left(  \cup_{t\in\tilde{T}_{\varepsilon/(2s)}(x),s>0\text{ }%
}\partial_{\varepsilon/2}(s\bar{f}_{t})(x)\right)  \cup\left(  \cup_{t\in
T\setminus\mathcal{P}\text{ }}\mathrm{N}_{\operatorname*{dom}\bar{f}_{t}%
}^{\varepsilon/2}(x)\right)  \right)  \right]  _{\infty},\label{arsal}%
\end{equation}
where $\tilde{T}_{\varepsilon/(2s)}(x):=\{t\in\mathcal{P}:s\bar{f}_{t}%
(x)\geq-\varepsilon/2\}.$ Moreover, since $x\in\operatorname*{dom}f$
$(\subset\operatorname*{dom}f_{t},$ for all $t\in T),$ we have
\[
\mathrm{N}_{\operatorname*{dom}f_{t}}^{\varepsilon/2}(x)=\mathrm{N}%
_{\overline{\operatorname*{dom}f_{t}}}^{\varepsilon/2}(x)=\mathrm{N}%
_{\overline{\operatorname*{dom}\bar{f}_{t}}}^{\varepsilon/2}(x)=\mathrm{N}%
_{\operatorname*{dom}\bar{f}_{t}}^{\varepsilon/2}(x),\text{ for all}\ t\in T.
\]
At the same time, since $x\in\lbrack f\leq0]=\cap_{t\in T}[f_{t}\leq0],$ for
each $t\in\tilde{T}_{\varepsilon/(2s)}(x)\,$we have
\[
sf_{t}(x)\geq s\bar{f}_{t}(x)\geq-\varepsilon/2\geq sf_{t}(x)-\varepsilon/2,
\]
which implies that $t\in T_{\varepsilon/(2s)}(x)\subset T_{\varepsilon/s}(x)$
and $\partial_{\varepsilon/2}(s\bar{f}_{t})(x)\subset\partial_{\varepsilon
}(sf_{t})(x).$ In other words,
\[
\cup_{t\in\tilde{T}_{\varepsilon/(2s)}(x),s>0\text{ }}\partial_{\varepsilon
/2}(s\bar{f}_{t})(x)\subset\cup_{t\in T_{\varepsilon/s}(x),s>0\text{ }%
}\partial_{\varepsilon}(sf_{t})(x),
\]
and (\ref{arsal}) entails the inclusion \textquotedblleft$\subset
$\textquotedblright\ in (\ref{a1a}). 

We proceed by showing the opposite inclusion \textquotedblleft$\supset
$\textquotedblright\ in (\ref{a1a}). Indeed, if $A_{\varepsilon}$ and
$B_{\varepsilon}$ are as in (\ref{yarab2}), then for every $x^{\ast}\in
A_{\varepsilon}\cup B_{\varepsilon}$ and $y\in\lbrack f\leq0]$ $(\subset
\operatorname*{dom}f_{t},$ for all $t\in T)$ we have\
\[
\left\langle x^{\ast},y-x\right\rangle \leq\max\{\sup_{t\in T_{\varepsilon
/s}(x),s>0\text{ }}(sf_{t}(y)-sf_{t}(x))+\varepsilon,\varepsilon
\}\leq2\varepsilon.
\]
The same inequality holds when\ $x^{\ast}\in\overline{\operatorname*{co}%
}\left(  A_{\varepsilon}\cup B_{\varepsilon}\right)  ,\ $and therefore\
\[
\sigma_{\left[  \overline{\operatorname*{co}}\left(  A_{\varepsilon}\cup
B_{\varepsilon}\right)  \right]  _{\infty}}(y-x)\leq0,\text{ for all }%
y\in\lbrack f\leq0].
\]
Hence, $\left[  \overline{\operatorname*{co}}\left(  A_{\varepsilon}\cup
B_{\varepsilon}\right)  \right]  _{\infty}\subset\mathrm{N}_{[f\leq0]}(x),$
and the proof of (\ref{a1a}) is complete.

Finally, to prove (\ref{a2b}), we observe that for every $\varepsilon>0$ and
$\delta>0\ $%
\[
\delta(A_{\varepsilon}\cup B_{\varepsilon})=\left(  \cup_{t\in T_{\varepsilon
/s}(x),s>0\text{ }}\partial_{\delta\varepsilon}(\delta sf_{t})(x)\right)
\cup\left(  \cup_{t\in T\setminus\mathcal{P}\text{ }}\mathrm{N}%
_{\operatorname*{dom}f_{t}}^{\delta\varepsilon}(x)\right)  =A_{\delta
\varepsilon}\cup B_{\delta\varepsilon},
\]
where $A_{\varepsilon}$ and $B_{\varepsilon}$ are defined as in (\ref{yarab2}%
). In other words, the set $\cap_{\varepsilon}\overline{\operatorname*{co}%
}\left(  A_{\varepsilon}\cup B_{\varepsilon}\right)  $ is a (convex) cone.
Therefore, invoking (\ref{ricone}), (\ref{a2b}) follows from (\ref{a1a}) by
taking the intersection over all $\varepsilon>0$.
\end{dem}

\begin{rem}
\label{remm}$(i)$ Related formulas to the characterization given in
(\ref{a2b}) were obtained\ in\ \cite[Corollary 7]{HaSv17} for the special case
where $f(x)=0$. In contrast, the formulas derived in \cite[Example
5.1.13]{CHLBook} are expressed\ in terms of convex combinations of the data
functions and explicitly involve the domain $\operatorname*{dom}f$. In both
references, the functions $f_{t}$ are assumed to be proper and lsc.

$(ii)$ Condition (\ref{cc2}) is clearly satisfied\ when\ all the $f_{t}$'s are
lsc, but it also holds in more general situations.\ For instance, if $f$ is
continuous at some$\ $point in $[f\leq0]$, then $\overline{f}=\sup_{t\in
T}\bar{f}_{t}$ by \cite[Proposition 5.2.4]{CHLBook}, and condition\ (\ref{cc2}%
) follows from (\ref{lemnewaa}). Moreover, the equality $\overline{[f\leq
0]}=[\overline{f}\leq0]$ also holds whenever\ $[f<0]\cap f^{-1}(\mathbb{R}%
)\neq\emptyset,\ $again by (\ref{lemnewaa}).
\end{rem}

\medskip In the particular case of a singleton family, the preceding theorem
takes the simpler form (see, also, \cite{HaSv17}). We refer to \cite{Ku79} and
references therein for the origins of these formulas.

\begin{cor}
Let\ $f:X\rightarrow\mathbb{R}_{\infty}$ be a convex$\ $function such that
$\bar{f}$ is proper and $\operatorname*{cl}([f\leq0])=[\operatorname*{cl}%
f\leq0]$. Let $x\in\lbrack f\leq0]$ be such that $f(x)\in\mathbb{R}.$ Then
\[
\mathrm{N}_{[f\leq0]}(x)=\left[  \operatorname*{cl}\left(
{\textstyle\bigcup\limits_{sf(x)\geq-\varepsilon,s>0\text{ }}}
\partial_{\varepsilon}(sf)(x)\right)  \right]  _{\infty},\text{ for all
}\varepsilon>0,
\]
and, consequently,
\[
\mathrm{N}_{[f\leq0]}(x)=%
{\textstyle\bigcap\nolimits_{\varepsilon>0}}
\operatorname*{cl}\left(
{\textstyle\bigcup\limits_{sf(x)\geq-\varepsilon,s>0\text{ }}}
\partial_{\varepsilon}(sf)(x)\right)  .
\]

\end{cor}

\begin{dem}
Since $\bar{f}$ is proper by assumption, Theorem \ref{thmcone} entails
\[
\mathrm{N}_{[f\leq0]}(x)=\left[  \overline{\operatorname*{co}}\left(
\cup_{sf(x)\geq-\varepsilon,s>0}\partial_{\varepsilon}(sf)(x)\right)  \right]
_{\infty}.
\]
Thus, since the set $\cup_{sf(x)\geq-\varepsilon,s>0}\partial_{\varepsilon
}(sf)(x)$ is convex, the first equality of the assertion holds. The second
equality is obtained as in the final\ part of the proof of Theorem
\ref{thmcone}.
\end{dem}

The following corollary provides a characterization of the normal cone to
$\operatorname*{cl}([f<0])$ at points belonging to $\operatorname*{cl}%
([f<0]),$ without requiring them to lie in $[f<0]$. Such cones are important
in optimization theory and variational analysis (see, e.g., \cite{Da}).

\begin{cor}
\label{cor10}Let\ $f_{t}:X\rightarrow\overline{\mathbb{R}}$, $t\in T,$ be
convex functions,$\ $and suppose\ that $f:=\sup_{t\in T}f_{t}$ is continuous
at some point in $[f<0]\cap f^{-1}(\mathbb{R}).$ Then, for every $x\in\lbrack
f\leq0]$ such that $f(x)\in\mathbb{R},$ we have
\[
\mathrm{N}_{\operatorname*{cl}([f<0])}(x)=\left[  \overline{\operatorname*{co}%
}\left(  \left(
{\textstyle\bigcup\nolimits_{t\in T_{\frac{\varepsilon}{s}}(x),s>0\text{ }}}
\partial_{\varepsilon}(sf_{t})(x)\right)  \cup\left(
{\textstyle\bigcup\nolimits_{t\in T\setminus\mathcal{P}\text{ }}}
\mathrm{N}_{\operatorname*{dom}f_{t}}^{\varepsilon}(x)\right)  \right)
\right]  _{\infty},
\]
for all $\varepsilon>0.$ Consequently,
\[
\mathrm{N}_{\operatorname*{cl}([f<0])}(x)=%
{\textstyle\bigcap\nolimits_{\varepsilon>0}}
\overline{\operatorname*{co}}\left(  \left(
{\textstyle\bigcup\nolimits_{t\in T_{\frac{\varepsilon}{s}}(x),s>0\text{ }}}
\partial_{\varepsilon}(sf_{t})(x)\right)  \cup\left(
{\textstyle\bigcup\nolimits_{t\in T\setminus\mathcal{P}\text{ }}}
\mathrm{N}_{\operatorname*{dom}f_{t}}^{\varepsilon}(x)\right)  \right)  .
\]

\end{cor}

\begin{dem}
Fix $x\in\lbrack f\leq0]$ with $f(x)\in\mathbb{R}.$ By \cite[Proposition
5.2.4]{CHLBook}, the continuity assumption on $f$ ensures that $\bar{f}%
=\sup_{t\in T}\bar{f}_{t}.$ Therefore, by (\ref{lemnewa}), the condition
$[f<0]\cap f^{-1}(\mathbb{R})\neq\emptyset$ entails\
\[
\operatorname*{cl}([f<0])=\operatorname*{cl}([f\leq0])=[\bar{f}\leq
0]=[\sup_{t\in T}\bar{f}_{t}\leq0]=\cap_{t\in T}[\bar{f}_{t}\leq0],
\]
whish shows that the family $\{f_{t},$ $t\in T\}$ satisfies condition
(\ref{cc2}). Consequently, since $x\in\lbrack f\leq0]$ we have
\[
\mathrm{N}_{\operatorname*{cl}([f<0])}(x)=\mathrm{N}_{\operatorname*{cl}%
([f\leq0])}(x)=\mathrm{N}_{[f\leq0]}(x),
\]
and the conclusion follows by applying Theorem \ref{thmcone}.
\end{dem}

We now apply the previous results on the normal cone to\ derive
general\ optimality conditions for the following convex infinite convex
program, posed in the lcs $X$,\
\[
(\mathtt{P})\text{\ \ }\inf_{f_{t}(x)\leq0,\text{ }t\in T}f_{0}(x),
\]
where each function $f_{t}:X\rightarrow\overline{\mathbb{R}}$, $t\in
T\cup\{0\}$ (with $T\neq\emptyset$ and $0\notin T$), is assumed to be convex.
We suppose that problem $(\mathtt{P})$ has a finite optimal value
$v(\mathtt{P})$, and we denote\
\[
f:=\sup_{t\in T}f_{t}.
\]
Unlike classical optimality conditions, which typically require additional
qualifications such as Slater's condition (see, i.e., \cite[Theorem
23.7]{Ro70}), our approach relies only on the lower semicontinuity-like
assumption (\ref{cc2}), namely,
\[
\overline{\lbrack f\leq0]}=\cap_{t\in T}[\overline{f_{t}}\leq0].
\]

The result below extends the optimality rules previously obtained in the
compact-continuous framework in \cite{HaLo23} (see, also, \cite{CompactCase}).

\begin{theo}
\label{thmopt} Let\ $f_{t}:X\rightarrow\overline{\mathbb{R}}$, $t\in
T\cup\{0\},$ be convex functions such that $v(\mathtt{P})\in\mathbb{R}$, and
denote $f:=\sup_{t\in T}f_{t}.$ Assume that either $f_{0}$ is continuous at
some point in $[f\leq0]$ or $\operatorname*{int}([f\leq0])\cap
\operatorname*{dom}f_{0}\neq\emptyset.$ If the family $\{f_{t},$ $t\in T\}$
satisfies condition (\ref{cc2}), then a feasible point $x$ of\ $(\mathtt{P})$
is optimal if and only if$,$ for every $\varepsilon>0,$
\[
\theta\in\partial f_{0}(x)+\left[  \overline{\operatorname*{co}}\left(
\left(
{\textstyle\bigcup\nolimits_{t\in T_{\frac{\varepsilon}{s}}(x),s>0\text{ }}}
\partial_{\varepsilon}(sf_{t})(x)\right)  \cup\left(
{\textstyle\bigcup\nolimits_{t\in T\setminus\mathcal{P}\text{ }}}
\mathrm{N}_{\operatorname*{dom}f_{t}}^{\varepsilon}(x)\right)  \right)
\right]  _{\infty},
\]
or, equivalently,
\[
\theta\in\partial f_{0}(x)+%
{\textstyle\bigcap\nolimits_{\varepsilon>0}}
\overline{\operatorname*{co}}\left(  \left(
{\textstyle\bigcup\nolimits_{t\in T_{\frac{\varepsilon}{s}}(x),s>0\text{ }}}
\partial_{\varepsilon}(sf_{t})(x)\right)  \cup\left(
{\textstyle\bigcup\nolimits_{t\in T\setminus\mathcal{P}\text{ }}}
\mathrm{N}_{\operatorname*{dom}f_{t}}^{\varepsilon}(x)\right)  \right)  .
\]

\end{theo}

\begin{dem}
A feasible point $x$ is optimal for\ $(\mathtt{P})$ if and only if, by the
Moreau-Rockafellar theorem (e.g., \cite[Ch. 3]{CHLBook}),
\[
\theta\in\partial(f_{0}+\mathrm{I}_{[f\leq0]})(x)=\partial f_{0}%
(x)+\mathrm{N}_{[f\leq0]}(x).
\]
The conclusion then follows directly\ from Theorem \ref{thmcone}.$\smallskip$
\end{dem}

\medskip We give an illustration of Theorem \ref{thmopt} in the setting of
infinite linear optimization, where $(\mathtt{P})$ takes the form
\begin{equation}
\inf_{\left\langle a_{t},x\right\rangle \leq b_{t},\text{ }t\in T}\left\langle
c,x\right\rangle ,\label{lin}%
\end{equation}
with\ $c,a_{t}\in X^{\ast},$ $t\in T,$ being given vectors. Applying Theorem
\ref{thmopt}, we obtain the following result.

\begin{cor}
Consider the infinite linear program (\ref{lin}). Then a feasible point\ $x\in
X$ is optimal if and only if
\begin{equation}
-c\in\left[  \overline{\operatorname*{co}}\{sa_{t}:\left\langle a_{t}%
,x\right\rangle -b_{t}\geq-\varepsilon/s,\text{ }t\in T,\text{ }s>0\}\right]
_{\infty},\text{ for every }\varepsilon>0,\label{ll}%
\end{equation}
or, equivalently,
\[
-c\in%
{\textstyle\bigcap\limits_{\varepsilon>0}}
\overline{\operatorname*{co}}\{sa_{t}:\left\langle a_{t},x\right\rangle
-b_{t}\geq-\varepsilon/s,t\in T,s>0\}.
\]

\end{cor}

If $T$ is finite in the previous corollary, then (\ref{ll}) reduces\ to the
classical condition
\[
-c\in\operatorname*{cone}\{a_{t}:\left\langle a_{t},x\right\rangle
=b_{t},\text{ }t\in T\}.
\]

\section{Normal cone to sublevel sets in the quasi-convex case\label{sec5}}

In this section, we extend the previous results to the quasi-convex setting,
providing new characterizations of the normal cone to sublevel sets and
deriving optimality conditions for infinite quasi-convex optimization
problems. 

We consider an arbitrary family of quasi-convex functions $f_{t}%
:X\rightarrow\overline{\mathbb{R}}$, $t\in T$, and define the associated
supremum function by
\[
f:=\sup_{t\in T}f_{t}.
\]
By means of Theorem \ref{thmcone}, we then obtain initial characterizations of
the normal cone to $[f\leq0]$ in terms of the approximate normal sets
\[
\mathrm{N}_{[f_{t}\leq0]}^{\varepsilon}(x),\text{ }t\in T.
\]
We first present a more geometric reformulation of the closedness condition
(\ref{cc2}). The resulting condition (\ref{cc3}) is satisfied under natural
assumptions within the quasi-convex framework (see Lemma \ref{lemnewb}). 

\begin{lem}
Let\ $f_{t}:X\rightarrow\overline{\mathbb{R}},$ $t\in T,$ be given functions,
and denote $f:=\sup_{t\in T}f_{t}.$ If condition (\ref{cc2}) holds, namely,
$\overline{[f\leq0]}=\cap_{t\in T}[\bar{f}_{t}\leq0],$ then
\begin{equation}
\overline{\lbrack f\leq0]}=\cap_{t\in T}\operatorname*{cl}[f_{t}%
\leq0].\label{cc3}%
\end{equation}

\end{lem}

\begin{dem}
For every\ $t\in T,$ we have\ $\operatorname*{cl}([f_{t}\leq0])\subset
\lbrack\overline{f}_{t}\leq0]$ and $[f\leq0]\subset\lbrack f_{t}\leq0].$
Hence, condition (\ref{cc2}) yields
\[
\cap_{t\in T}\operatorname*{cl}([f_{t}\leq0])\subset\cap_{t\in T}[\overline
{f}_{t}\leq0]=\overline{[f\leq0]}\subset\cap_{t\in T}\operatorname*{cl}%
([f_{t}\leq0]).
\]
In other words, we have
\[
\overline{\lbrack f\leq0]}=\cap_{t\in T}\operatorname*{cl}([f_{t}\leq
0])=\cap_{t\in T}[\overline{f}_{t}\leq0],
\]
which\ proves (\ref{cc3}).
\end{dem}

We give\ the first characterizations of $\mathrm{N}_{[f\leq0]}(x)$ in the
quasi-convex setting.

\begin{theo}
\label{qthmcone} Let\ $f_{t}:X\rightarrow\overline{\mathbb{R}}$, $t\in T,$ be
quasi-convex functions, and denote $f:=\sup_{t\in T}f_{t}.$ Assume that
condition (\ref{cc3}) holds. Then, for every $x\in\lbrack f\leq0]$ such that
$f(x)\in\mathbb{R},$ we have\
\[
\mathrm{N}_{[f\leq0]}(x)=\left[  \overline{\operatorname*{co}}\left(
{\textstyle\bigcup\limits_{t\in T}}
\mathrm{N}_{[f_{t}\leq0]}^{\varepsilon}(x)\right)  \right]  _{\infty},\text{
for all }\varepsilon>0,
\]
and, consequently,
\[
\mathrm{N}_{[f\leq0]}(x)=%
{\textstyle\bigcap\limits_{\varepsilon>0}}
\overline{\operatorname*{co}}\left(
{\textstyle\bigcup\limits_{t\in T}}
\mathrm{N}_{[f_{t}\leq0]}^{\varepsilon}(x)\right)  .
\]

\end{theo}

\begin{dem}
Fix\ $x\in\lbrack f\leq0]\cap f^{-1}(\mathbb{R})$ and $\varepsilon>0.$
We\ introduce\ the proper convex functions defined on $X$ by
\[
g_{t}:=\mathrm{I}_{[f_{t}\leq0]},\text{ }t\in T,\text{ and }g:=\sup_{t\in
T}g_{t}.
\]
Then we have $g=\mathrm{I}_{\cap_{t\in T}[f_{t}\leq0]}=\mathrm{I}_{[f\leq0]}.$
Moreover, by (\ref{cc3}), we also have
\[
\bar{g}=\mathrm{I}_{\operatorname*{cl}([f\leq0])}=\sup_{t\in T}\mathrm{I}%
_{\operatorname*{cl}([f_{t}\leq0])}=\sup_{t\in T}\bar{g}_{t}.
\]
Consequently,
\[
\overline{\lbrack g\leq0]}=\operatorname*{cl}([f\leq0])=\cap_{t\in
T}\operatorname*{cl}([f_{t}\leq0])=\cap_{t\in T}[\bar{g}_{t}\leq0],
\]
and the family $\{g_{t},$ $t\in T\}$ satisfies condition (\ref{cc2}). Next,
observe that each function $\bar{g}_{t}$ is proper and satisfies $sg_{t}%
=g_{t},$ for all $s>0.$ Furthermore, since $x\in\lbrack f\leq0]$
$(\subset\lbrack f_{t}\leq0],$ for all $t\in T),$ it follows that
$g_{t}(x)=0,$ for all $t\in T.$ Thus, according to Theorem \ref{thmcone}, for
each $\varepsilon>0$ we have
\[
\mathrm{N}_{[f\leq0]}(x)=\mathrm{N}_{[g\leq0]}(x)=\left[  \overline
{\operatorname*{co}}\left(
{\textstyle\bigcup\limits_{t\in T}}
\partial_{\varepsilon}g_{t}(x)\right)  \right]  _{\infty}=\left[
\overline{\operatorname*{co}}\left(
{\textstyle\bigcup\limits_{t\in T}}
\mathrm{N}_{[f_{t}\leq0]}^{\varepsilon}(x)\right)  \right]  _{\infty},
\]
which establishes the first statement. The second assertion follows
analogously from Theorem \ref{thmcone}.
\end{dem}

The following corollary provides sufficient conditions\ for the validity of
the closure condition\ (\ref{cc3}), thus enabling the application of Theorem
\ref{qthmcone}.

\begin{cor}
Let\ $f_{t}:X\rightarrow\overline{\mathbb{R}}$, $t\in T,$ be quasi-convex
functions, and denote $f:=\sup_{t\in T}f_{t}.$ Assume that $[f\leq
0]\neq\emptyset.$ Then the conclusion of Theorem \ref{qthmcone} remains valid
under either of the following conditions:

$(i)$ Function $f$ is continuous at some point in $[f\leq0].$

$(ii)$ $X=\mathbb{R}^{n},$ $\bar{f}=\sup_{t\in T}\bar{f}_{t}$,
$\operatorname*{ri}(\operatorname*{dom}f)\cap\lbrack f\leq0]\neq\emptyset,$
and $\operatorname*{ri}(\operatorname*{dom}f_{t})\cap\lbrack f\leq
0]\neq\emptyset$ for all $t\in T.$
\end{cor}

\begin{dem}
Assume that condition $(i)$ holds. Since all the functions $\mathrm{I}%
_{[f\leq0]}$ and $\mathrm{I}_{[f_{t}\leq0]},$ $t\in T,$ are convex, and
$\mathrm{I}_{[f\leq0]}=\sup_{t\in T}\mathrm{I}_{[f_{t}\leq0]},$ the continuity
assumption implies that (\cite[Proposition 5.2.4]{CHLBook})
\[
\mathrm{I}_{\operatorname*{cl}([f\leq0])}=\operatorname*{cl}(\mathrm{I}%
_{[f\leq0]})=\sup_{t\in T}\operatorname*{cl}(\mathrm{I}_{[f_{t}\leq0]}%
)=\sup_{t\in T}\mathrm{I}_{\operatorname*{cl}([f_{t}\leq0])}.
\]
Hence, $\operatorname*{cl}([f\leq0])=\cap_{t\in T}\operatorname*{cl}%
([f_{t}\leq0])$ and condition (\ref{cc3}) holds. Thus, the conclusion follows
from Theorem \ref{qthmcone}.

Under condition $(ii),$ the relation\ $\bar{f}=\sup_{t\in T}\bar{f}_{t}$
implies that
\[
\lbrack\overline{f}\leq0]=\cap_{t\in T}[\overline{f}_{t}\leq0].
\]
Moreover, since $\emptyset\neq\lbrack f\leq0]\subset\lbrack f_{t}\leq0],$ for
all $t\in T,$ Lemma \ref{lemnewb} entails $[\overline{f}\leq
0]=\operatorname*{cl}[f\leq0]$ and $[\overline{f}_{t}\leq0]=\operatorname*{cl}%
[f_{t}\leq0],$ for all $t\in T.$ Thus, combining these identities yields
condition (\ref{cc3}). The conclusion of the corollary then follows directly
from Theorem \ref{qthmcone}.
\end{dem}

To provide a more explicit description of the normal sets $\mathrm{N}%
_{[f_{t}\leq0]}^{\varepsilon}(x)$ in Theorem \ref{qthmcone}, we give the
following lemma in the context of Asplund spaces--i.e., Banach spaces in which
the dual of every separable subspace is separable. This class includes all
reflexive Banach spaces. 

Here, $B_{X}(x,r)$ denotes the closed ball in $X$ centered at $x$ with radius
$r>0,$ and $\operatorname*{cl}\nolimits^{w^{\ast}}$ refers to the closure in
$X^{\ast}$ with respect to the weak*-topology.

\begin{lem}
\label{lemnos}Assume that $X$ is an Asplund Banach space. Let $f:X\rightarrow
\mathbb{R}_{\infty}$ be a lsc quasi-convex function. Then, for every
$x\in\lbrack f\leq0]$ and $\varepsilon>0,$ we have
\[
\mathrm{N}_{[f\leq0]}^{\varepsilon}(x)\subset\operatorname*{cl}%
\nolimits^{w^{\ast}}\{\mathbb{R}_{+}f(y):y\in B_{X}(x,3\sqrt{\varepsilon
}),\text{ }f(y)\leq2\sqrt{\varepsilon}\}+B_{X^{\ast}}(\theta,\sqrt
{\varepsilon}).
\]

\end{lem}

\begin{dem}
Fix $x\in\lbrack f\leq0]$ and pick $x^{\ast}\in\mathrm{N}_{[f\leq
0]}^{\varepsilon}(x),$ for $\varepsilon>0$. Since the indicator function
$\mathrm{I}_{[f\leq0]}$ is convex, proper, and lsc, the Brøndsted-Rockafellar
theorem (e.g., \cite[Proposition 4.3.8]{CHLBook}) guarantees the existence of
some $x_{\varepsilon}\in B_{X}(x,\sqrt{\varepsilon})$ and $x_{\varepsilon
}^{\ast}\in B_{X^{\ast}}(x^{\ast},\sqrt{\varepsilon})\cap\partial
(\mathrm{I}_{[f\leq0]})(x_{\varepsilon})$ such that $\left\vert \left\langle
x^{\ast},x_{\varepsilon}-x\right\rangle \right\vert \leq2\varepsilon.$
Moreover, since
\[
\mathrm{I}_{[f\leq0]}=\mathrm{I}_{[f^{+}\leq0]}=\sup_{k\geq1}kf^{+},
\]
and the sequence $(kf^{+})_{k\geq1}$ is a non-decreasing sequence of lsc
functions, by \cite[Proposition 3.7]{PP} it follows that
\[
x_{\varepsilon}^{\ast}\in\operatorname*{cl}\nolimits^{w^{\ast}}\{k\partial
f^{+}(u):u\in B_{X}(x_{\varepsilon},\sqrt{\varepsilon}),\text{ }f^{+}(u)\leq
k^{-1}\sqrt{\varepsilon}\}.
\]
For each such $u$, we also have (see \cite[Proposition 3.2(v)]{PP})
\[
\partial f^{+}(u)\subset\operatorname*{cl}\nolimits^{w^{\ast}}\{[0,1]\partial
f(u_{1}):u_{1}\in B_{X}(u,\sqrt{\varepsilon}),\text{ }\left\vert
f(u_{1})-f(u)\right\vert \leq\sqrt{\varepsilon}\}.
\]
Consequently, we obtain
\[
x_{\varepsilon}^{\ast}\in\operatorname*{cl}\nolimits^{w^{\ast}}\{\lambda
k\partial f(u_{1}):u_{1}\in B_{X}(x_{\varepsilon},2\sqrt{\varepsilon}),\text{
}f(u_{1})\leq(1+k^{-1})\sqrt{\varepsilon},\text{ }\lambda\in\lbrack0,1],\text{
}k\geq1\},
\]
and the desired statement follows as
\[
x^{\ast}\in\operatorname*{cl}\nolimits^{w^{\ast}}\{\mathbb{R}_{+}f(y):y\in
B_{X}(x,3\sqrt{\varepsilon}),\text{ }f(y)\leq2\sqrt{\varepsilon}\}+B_{X^{\ast
}}(\theta,\sqrt{\varepsilon}).
\]

\end{dem}

We give the second main result of this section. 

\begin{theo}
\label{qthmconeb} Assume that $X$ is an Asplund space. Let\ $f_{t}%
:X\rightarrow\overline{\mathbb{R}}$, $t\in T,$ be lsc quasi-convex functions,
and denote $f:=\sup_{t\in T}f_{t}.$ Then, for every $x\in\lbrack f\leq0]$ such
that $f(x)\in\mathbb{R},$ we have\
\[
\mathrm{N}_{[f\leq0]}(x)\subset%
{\textstyle\bigcap\limits_{\varepsilon>0}}
\overline{\operatorname*{co}}\{\mathbb{R}_{+}\partial f_{t}(y):y\in
B_{X}(x,\varepsilon),\text{ }f_{t}(y)\leq\varepsilon,\text{ }t\in T\}.
\]

\end{theo}

\begin{dem}
Fix\ $x\in\lbrack f\leq0]\cap f^{-1}(\mathbb{R})$ and small enough
$\varepsilon>0.$ Since the functions $f_{t},$ $t\in T,$ are lsc, Theorem
\ref{qthmcone} applies and yields
\[
\mathrm{N}_{[f\leq0]}(x)=%
{\textstyle\bigcap\limits_{\varepsilon>0}}
\overline{\operatorname*{co}}\left(
{\textstyle\bigcup\limits_{t\in T}}
\mathrm{N}_{[f_{t}\leq0]}^{\varepsilon}(x)\right)  .
\]
Therefore, by Lemma \ref{lemnos},
\begin{align*}
\mathrm{N}_{[f\leq0]}(x)  & \subset%
{\textstyle\bigcap\limits_{\varepsilon>0}}
\overline{\operatorname*{co}}\left(
{\textstyle\bigcup\limits_{t\in T}}
\operatorname*{cl}\nolimits^{w^{\ast}}\{\mathbb{R}_{+}f_{t}(y):y\in
B_{X}(x,3\sqrt{\varepsilon}),\text{ }f_{t}(y)\leq2\sqrt{\varepsilon
}\}+B_{X^{\ast}}(\theta,\sqrt{\varepsilon})\right)  \\
& =%
{\textstyle\bigcap\limits_{\varepsilon>0}}
\overline{\operatorname*{co}}\left(
{\textstyle\bigcup\limits_{t\in T}}
\{\mathbb{R}_{+}\partial f_{t}(y):y\in B_{X}(x,\varepsilon),\text{ }%
f_{t}(y)\leq\varepsilon\}\right)  ,
\end{align*}
and the result follows. 
\end{dem}

Next, we consider the infinite optimization problem $(\mathtt{P})$ in the
Asplund space $X:$\
\[
(\mathtt{P})\text{\ \ }\inf_{f_{t}(x)\leq0,\text{ }t\in T}f_{0}(x),
\]
where each $f_{t}:X\rightarrow\mathbb{R}_{\infty}$, $t\in T$, is lsc and
quasi-convex. We suppose that problem $(\mathtt{P})$ has a finite optimal
value $v(\mathtt{P})$, and we denote\
\[
f:=\sup_{t\in T}f_{t}.
\]
We assume that $f_{0}$ is convex and that either $f_{0}$ is continuous at some
point in $[f\leq0]$ or $\operatorname*{int}([f\leq0])\cap\operatorname*{dom}%
f_{0}\neq\emptyset.$ The following theorem establishes a necessary optimality
condition for problem $(\mathtt{P})$ in this quasi-convex setting.

\begin{theo}
\label{thmoptbb} If $x$ is feasible for\ $(\mathtt{P}),$ then we have
\[
\theta\in\partial f_{0}(x)+%
{\textstyle\bigcap\limits_{\varepsilon>0}}
\overline{\operatorname*{co}}\{\mathbb{R}_{+}\partial f_{t}(y):y\in
B_{X}(x,\varepsilon),\text{ }f_{t}(y)\leq\varepsilon,\text{ }t\in T\}.
\]

\end{theo}

\begin{dem}
We have
\[
\theta\in\partial(f_{0}+\mathrm{I}_{[f\leq0]})(x)=\partial f_{0}%
(x)+\mathrm{N}_{[f\leq0]}(x),
\]
and the conclusion follows by Theorem \ref{qthmconeb}.$\smallskip$
\end{dem}

\section{Conclusion\label{sec6}}

We have developed new characterizations of the normal cone to the sublevel set
of the supremum of an arbitrary family of convex and quasi-convex functions.
In the convex setting, these results are expressed exclusively in terms of
subdifferentials of the individual data functions at the nominal point. In
contrast, the quasi-convex framework requires the use of the Fréchet
subdifferential of these individual data functions but evaluated at nearby points.\ 

These results are applied to derive necessary optimality conditions for
infinite convex and quasi-convex optimization problems.

Our approach relies on techniques from convex and nonsmooth analysis. The
quasi-convexity assumption allows the reduction of the problem---whether for
subdifferential calculus or for deriving optimality conditions---to the
description of the normal cone to a sublevel set of a quasi-convex function.

\backmatter

%
%


\bmhead{Acknowledgements}

The research is supported by the VRII-UNAP Internal Project
2024-2025, Grant PID2022-136399NB-C21 funded by
MICIU/AEI/10.13039/501100011033 and by ERDF/EU, and\ Basal CMM FB210005. 

\section*{Declarations}


\begin{itemize}
\item Conflict of interest/Competing interests. Not applicable. 
\item Ethics approval and consent to participate. Not applicable.
\item Data availability. Not applicable. 
\end{itemize}


\begin{thebibliography}{99}                                                                                               %
\bibitem {CaTh14}\textsc{Cabot, A., Thibault, L., }Sequential formulae for the
normal cone to sublevel sets. Trans. Amer. Math. Soc. 366 (2014) 6591--6628.

\bibitem {CaHa26}\textsc{Caro, S., Hantoute, A.,} Characterization of the
normal cone to the domain and optimality condition in infinite optimization. (2026).

\bibitem {CompactCase}\textsc{Correa, R., Hantoute, A., López, M. A.,}
Alternative representations of the normal cone to the domain of supremum
functions and subdifferential calculus. Set-Valued Var. Anal. 29 (2021) 683--699.

\bibitem {CHLBook}\textsc{Correa, R., Hantoute, A., López, M. A.,}
\emph{Fundamentals of Convex Analysis and Optimization, A Supremum Function
Approach}, Springer Cham, (2023).

\bibitem {Cr2005}\textsc{Crouzeix, J. P.,} Continuity and Differentiability of
Quasiconvex Functions. In: Hadjisavvas, N., Koml si, S., Schaible, S. (eds)
Handbook of Generalized Convexity and Generalized Monotonicity. Nonconvex
Optimization and Its Applications, vol 76. Springer, New York, NY, 2005.

\bibitem {Da}\textsc{Daniilidis, A., Hadjisavvas, N., Martínez-Legaz, J. E.,}
An Appropriate Subdifferential for Quasiconvex Functions. SIAM J. Optim. 12
(2001) 407--420.

\bibitem {DiGoLo06}\textsc{Dinh, N., Goberna, M.A., López, M.A., } From linear
to convex systems: consistency, Farkas' lemma and applications, J. Convex
Anal. 13 (2006) 113--133.

\bibitem {HaLo23}Hantoute, A., López, M.A., New tour on the subdifferential of
supremum via finite sums and suprema. J. Optim. Theory App. 193 (2022) 167--194.

\bibitem {HLZ08}\textsc{Hantoute, A., López, M. A., Z\u{a}linescu,
C.},\ Subdifferential calculus rules in convex analysis: A unifying approach
via pointwise supremum functions. SIAM J. Optim. 19 (2008) 863--882.

\bibitem {HaSv17}\textsc{Hantoute, A., Svensson, A.,} A General representation
of $\delta$-normal sets to sublevels of convex functions. Set-Valued Var. Anal
25 (2017) 651--678.

\bibitem {HiMoVoSe39}\textsc{Hiriart-Urruty, J.-B., Moussaoui, M., Seeger, A.,
Volle, M.}, Subdifferential calculus without qualification conditions, using
approximate subdifferentials: a survey. Nonlinear Anal. 24 (1995) 1727--1754.

\bibitem {Ku79}\textsc{Kutateladze S. S.}, \emph{Convex E-programming}. Soviet
Math. Dokl. 20 (1979) 391--393.

\bibitem {Mor}\textsc{Mordukhovich B. S.}, \emph{Variational Analysis and
Applications}. Springer (2018).

\bibitem {PP}\textsc{Pérez-Aros, P., }Subdifferential formulae for the
supremum of an arbitrary family of functions. SIAM J. Optim. 29 (2019) 1714--1743.

\bibitem {Ro70}\textsc{Rockafellar, R. T.} \emph{Convex Analysis}. Princeton
Mathematical Series Princeton University Press, Princeton, N. J. (1970).

\bibitem {RoWe}\textsc{Rockafellar, R.T., Wets, R. J.-B.}, \emph{Variational
Analysis} (Grundlehren der mathematischen Wissenschaften), Springer, (2010).
\end{thebibliography}
\end{document}